\documentclass[12pt,a4paper]{amsart}
\usepackage{graphicx,amssymb}
\usepackage{amsmath,amssymb,amscd}
\input xy
\xyoption{all}
\usepackage[all]{xy}
\usepackage{tikz}

\usepackage{hyperref}

\newcommand{\lra}{\longrightarrow}

\newcommand{\cO}{\mathcal{O}}

\newcommand{\cM}{\mathcal{M}}

\newcommand{\cL}{\mathcal L}
\newcommand{\cR}{\mathcal R}

\newcommand{\End}{\mbox{End}}

\newcommand{\Cliff}{\operatorname{Cliff}}
\newcommand{\gr}{\text{gr}}
\newcommand{\rk}{\text{rk}}

\newcommand{\Gr}{\operatorname{Gr}}
\newcommand{\go}{\operatorname{gon}}
\newcommand{\GL}{\operatorname{GL}}
\newcommand{\um}{\"}
\theoremstyle{plain}
\newtheorem{theorem}{Theorem}[section]

\newtheorem{conj}{Conjecture}
\newtheorem{rem}{Remark}[section]

\newtheorem{qu}[theorem]{Question}
\numberwithin{equation}{section}
\begin{document}
\title[Higher rank Brill-Noether theory]{Higher rank Brill-Noether theory and coherent systems\\Open Questions}

\author{P. E. Newstead}
\address{P. E. Newstead\\Department of Mathematical Sciences\\
              University of Liverpool\\
              Peach Street, Liverpool L69 7ZL, UK}
\email{newstead@liv.ac.uk}      
\date{\today}

\thanks{The author is a member of the research group VBAC (Vector Bundles on Algebraic Curves).}
\maketitle
\begin{abstract}
This article presents a list of open questions on higher rank Brill-Noether theory and coherent systems. Background material and appropriate references are included.
\end{abstract}
\tableofcontents

In this article, we present some of the many open questions on higher rank Brill-Noether theory and coherent systems on algebraic curves. A previous problems list including open questions on higher rank Brill-Noether theory was produced originally in 1991 with an updated version in 1994 for the foundational meeting of the research group Vector Bundles on Algebraic Curves (VBAC). Mercat's Pr\'esentation of 2001 contains a further list. In 2003, I launched the Brill-Noether project aimed at solving the basic questions in 10 years. Of course, this was not achieved (I never really expected it to be), but much progress has been made. The project and Mercat's Pr\'esentation are available on my website http://www.liverpool.ac.uk/$\sim$newstead; almost all the questions  stated there are still open in full generality, although they have been solved in specific cases. The survey article \cite{n1} concludes with a list of problems on moduli spaces of coherent systems and many open questions are stated in various papers on the subject.

The questions presented here are a mixture of general problems and more specific ones. They are far from being a complete list. The list of references is also incomplete; I have concentrated on those which form the basis for the questions raised in this article. Most of the questions are unlikely to have a clear-cut answer; partial solutions and counter-examples are a more plausible aim.

My thanks are due to Usha Bhosle, Leticia Brambila-Paz, Gavril Farkas, George Hitching, Michael Hoff, Emanuele Macr\'{\i}, Alexander Schmitt and Montserrat Teixidor i Bigas for comments, corrections, suggestions and additional references. Responsibility for any remaining errors and omissions is entirely mine.

\section{Background and definitions}\label{bd}

Brill-Noether theory is concerned with the possible dimensions of the space of sections $H^0(E)$ of a vector bundle on an algebraic variety. This is an extension of the classical study of meromorphic functions on a smooth complex algebraic curve (or Riemann surface) $C$ with poles along a given divisor. This study began with Riemann, Roch and Clifford, followed by a considerable number of other researchers in the 19th century, most notably for our purposes Brill and Max Noether. The essence of the theory is the study of projective embeddings of $C$. In more recent times, the theory has played a significant role in the study of moduli spaces of curves. 

We are concerned here with the higher-rank analogue; in other words, we consider vector bundles over a projective curve $C$, which we shall suppose to be defined over the complex numbers; many questions arise for curves defined over other fields (either not algebraically closed or of finite characteristic), but we shall not address these here. When $h^0(E)>\rk(E)$, there is a close connection with morphisms to Grassmannians and with Quot schemes, also with syzygies and with the geometrical structure of the moduli spaces of curves. In particular, counter-examples to the Harris-Morrison slope conjecture have been discovered using divisors in the moduli space which are defined in terms of rank-$2$ Brill-Noether loci (see, inter alia, \cite{fp} and \cite{fjp}).

In what follows, I have attempted to use a consistent set of notations. This is not the case in the literature.

Except in the final section, we shall suppose that $C$ is  smooth and irreducible. For $E$ a vector bundle, we write $\mu(E):=\frac{\deg E}{\rk E}$ for the \textit{slope} of $E$. Recall that a vector bundle $E$ is said to be \textit{stable} (\textit{semistable}) if every proper subbundle $F$ of $E$ satisfies $\mu(F)<(\le)\mu(E)$. Let $M(n,d)$ denote the moduli space of stable bundles of rank $n$ and degree $d$ on $C$, and $\widetilde{M}(n,d)$ the compactification of $M(n,d)$ whose points correspond to (S-equivalence classes of) semistable vector bundles; if $\gcd(n,d)=1$, $M(n,d)=\widetilde{M}(n,d)$. The \textit{Brill-Noether loci} $B(n,d,k)\subset M(n,d)$ and $\widetilde{B}(n,d,k)\subset\widetilde{M}(n,d)$ are defined by
\[B(n,d,k):=\{E\in M(n,d)|h^0(E)\ge k\}\]
and 
\[\widetilde{B}(n,d,k):=\{[E]\in \widetilde{M}(n,d)|h^0(\gr E)\ge k\},\]
where $[E]$ denotes the S-equivalence class of $E$ and $\gr E$ is the associated graded bundle. (In classical notation, $B(1,d,k)=W^{k-1}_d$; in the general situation, various other notations have been used.)

A \textit{coherent system} on $C$ \textit{of type} $(n,d,k)$ is a pair $(E,V)$, where $E$ is a vector bundle of rank $n$ and degree $d$, and $V$ is a linear subspace of $H^0(E)$ of dimension $k$; when $n=1$, these are referred to as \textit{linear systems}. When $k>n$, we say that $(E,V)$ is \textit{generated} (\textit{generically generated}) if the evaluation map $V\otimes \cO\to E$ is surjective (surjective at all but a finite number of points); if $(E,V)$ is generated, it determines a map from $C$ to the Grassmannian $\Gr(k,n)$. If $(E,H^0(E))$ is (generically) generated, we say also that $E$ is (generically) generated. 

For any $\alpha\in {\mathbb R}$, the $\alpha$-slope of $(E,V)$ is defined by
\[\mu_\alpha(E,V):=\frac{d+\alpha k}n.\]
A coherent system $(E,V)$ is said to be $\alpha$-\textit{stable} ($\alpha$-\textit{semistable}) if every proper subsystem $(F,W)$ of $(E,V)$ satisfies $\mu_\alpha(F,W)<(\le)\mu_\alpha(E,V)$. Note that $(E,0)$ is $\alpha$-stable ($\alpha$-semistable) if and only if $E$ is stable (semistable). It is easy to see that, if $k\ge1$, $(E,V)$ can be $\alpha$-semistable only if $\alpha\ge 0$, $d\ge0$ and $\alpha(n-k)\le d$ (note that the last of these conditions is vacuous if $k\ge n$). If $n=1$ and $\alpha>(\ge)0$, all coherent systems are $\alpha$-stable ($\alpha$-semistable). For any $\alpha$, there exists a quasi-projective moduli space $G(\alpha;n,d,k)$ of $\alpha$-stable coherent systems, which possesses a natural compactification $\widetilde{G}(\alpha;n,d,k)$ given by S-equivalence classes of $\alpha$-semistable coherent systems. When $n=1$, these spaces are independent of $\alpha>0$ and we write them as $G(1,d,k)$ (classically $G^{k-1}_d$).

A key feature for coherent systems is that of a \textit{critical value}. This can be defined as a value of $\alpha$ for which there exist $\alpha$-semistable coherent systems which become unstable (i.e. not semistable) for values of the parameter on one side or other  (possibly both) of $\alpha$. There is a simple numerical criterion for critical values, which are rational numbers with a finite set of possible denominators. It follows at once that there are finitely many critical values when $k<n$ and this can be proved also for $k\ge n$. We therefore have a sequence of critical values
\[0=\alpha_0<\alpha_1<...<\alpha_L<
\begin{cases}\frac{d}{n-k}& \text{if $k<n$},\\ \infty &\text{if $k\ge n$}.
\end{cases}\]
The moduli space can change only at a critical value and we write $G_i(n,d,k)$ ($G_L(n,d,k)$) for $G(\alpha;n,d,k)$ with $\alpha_i<\alpha<\alpha_{i+1}$ ($\alpha>\alpha_L$). The moduli space $G_L(n,d,k)$ may be referred to as the \textit{terminal} moduli space. One can define similarly $\widetilde{G}_i(n,d,k)$ and $\widetilde{G}_L(n,d,k)$. Further useful definitions are 
\[U^s(n,d,k):=\left\{(E,V)\in G_0(n,d,k)\left|\begin{array}{c}(E,V)\in G(\alpha;n,d,k) \text{ for }\\ \alpha>0, \alpha(n-k)<d\end{array}\right.\right\}\]
and $U(n,d,k):=\{(E,V)\in U^s(n,d,k)|E\text{ stable}\}$.

For any critical value $\alpha_i$, one can define \textit{flip loci} 
\[G_i^-:=\{(E,V)\in G_{i-1}(n,d,k)|(E,V)\not\in G_i(n,d,k)\})\] and 
\[G_i^+:=\{(E,V)\in G_i(n,d,k)|(E,V)\not\in G_{i-1}(n,d,k)\}.\]
These describe the wall-crossing at $\alpha_i$.

The \textit{Brill-Noether number} $\beta(n,d,k)$ is defined by 
\[\beta(n,d,k):=n^2(g-1)+1-k(k-d+n(g-1)),\]
where $g$ is the genus of $C$ (suppressed in the notation).
When $k\le d-n(g-1)$, $B(n,d,k)=M(n,d)$. Otherwise, the Brill-Noether number is often referred to as the \textit{expected dimension} of $B(n,d,k)$ and every irreducible component of $B(n,d,k)$ has dimension at least $\beta(n,d,k)$. Moreover, with no condition on $k$, every irreducible component of $G(\alpha;n,d,k)$ has dimension at least $\beta(n,d,k)$.

The infinitesimal behaviour of $G(\alpha;n,d,k)$ at $(E,V)$ is governed by the multiplication map
\[\mu_{E,V}:V\otimes H^0(E^*\otimes K)\lra H^0(E\otimes E^*\otimes K),\]
where $E^*$ is the dual of $E$ and $K$ is the canonical line bundle on $C$. This map is often referred to as the \textit{Petri map} at $(E,V)$. When $V=H^0(E)$, we write also $\mu_E$. The moduli space $G(\alpha;n,d,k)$ is smooth of dimension $\beta(n,d,k)$ at $(E,V)$ if and only if $\mu_{E,V}$ is injective. Moreover, if $E\in B(n,d,k)\setminus B(n,d,k+1)$, then $B(n,d,k)$ is smooth of dimension $\beta(n,d,k)$ at $E$ if and only if $\mu_E$ is injective. If $k>d-n(g-1)$, then $B(n,d,k+1)$ is contained in the singular set of $B(n,d,k)$.

If we consider bundles of fixed determinant $\cL$ of degree $d$, we have a moduli space $M(n,\cL)$, Brill-Noether loci $B(n,\cL,k)$ and moduli spaces $G(\alpha;n,\cL,k)$. It follows at once from the above that every irreducible component of $G(\alpha;n,\cL,k)$ has dimension at least $\beta(n,d,k)-g$, as does $B(n,\cL,k)$ when $k>d-n(g-1)$. However, statements about smoothness are not necessarily true. In particular, $B(2,K,k)$ has a new Brill-Noether number 
\[\beta(2,K,k):=3g-3-\frac{k(k+1)}2,\]
and a new Petri map 
\[\mu_E^K: S^2H^0(E)\lra H^0(S^2E),\]
arising from the symmetry $E\cong E^*\otimes K$.
All components of $B(2,K,k)$ have dimension at least $\beta(2,K,k)$ and, if $E\in B(2,K,k)\setminus B(2,K,k+1)$, then $B(2,K,k)$ is smooth of this dimension at $E$ if and only if $\mu_E^K$ is injective.

Given a vector bundle $F$ of rank $n'$ and degree $d'$, one can define the \textit{twisted Brill-Noether loci} $B(n,d,k)(F)$ by 
\[B(n,d,k)(F):=\{E\in M(n,d)|h^0(E\otimes F)\ge k\}\]
and similarly for $\widetilde{B}(n,d,k)(F)$. There is a new Brill-Noether number $\beta(n,d,k)(F)$ defined by
\[\beta(n,d,k)(F):=n^2(g-1)+1-k(k-n'd-nd'+nn'(g-1)) \]
and a new Petri map, defined as the composite of a multiplication map and a trace map:
\begin{eqnarray*}
\mu_{E_F}:H^0(F\otimes E)\otimes H^0(E^*\otimes F^*\otimes K)&\lra H^0(\End(F\otimes E)\otimes K)\\&\lra H^0(\End E\otimes K)
\end{eqnarray*}
with the same significance as before.
The simplest examples of such loci are generalised theta divisors, i.e. divisors of the form
\[\{E\in M(n,d)| h^0(E\otimes F)>0\},\]
where $F$ is a fixed vector bundle of rank $n'$ such that $\deg(E\otimes F)=nn'(g-1)$. These have been studied extensively over many years, but the extension of this to arbitrary twisted Brill-Noether loci is more recent. Here we could also allow $V$ to vary in a family, but we shall not pursue this further here.

One can extend this idea still further by choosing a homogeneous representation $\rho:\GL_n({\mathbb C})\to \GL(H)$, where $H$ is a vector space (this includes, for example, tensor power representations) and defining
\[B(\rho,n,d,k)(F):=\{E\in M(n,d)|h^0(E_\rho\otimes F)\ge k),\]
where $E_\rho$ is the vector bundle induced from $E$ by $\rho$. Again there is a similar definition for $\widetilde{B}(\rho,n,d,k)(F)$.
For all these twisted Brill-Noether loci, there is an obvious concept of coherent system.

\section{Basic questions}\label{bq}

Let $C$ be a smooth projective curve defined over $\mathbb C$. The following are basic questions in the Brill-Noether theory of $C$.

\begin {qu} For what values of $(n,d,k)$ is $B(n,d,k)$ non-empty?
\end{qu}
\begin {qu} What are the connected components of $B(n,d,k)$? In particular, for what values of $(n,d,k)$ is $B(n,d,k)$ connected?
\end{qu}
\begin{qu} What are the irreducible components of $B(n,d,k)$? In particular, for what values of $(n,d,k)$ is $B(n,d,k)$ irreducible? 
\end{qu}
\begin{qu} What is the dimension (of each component) of $B(n,d,k)$?
\end{qu}
\begin{qu} Does $B(n,d,k)$ have singularities outside $B(n,d,k+1)$? Could $B(n,d,k)$ be non-reduced?
\end{qu}

Similar questions can be asked for $\widetilde{B}(n,d,k)$, $B(n,\cL,k)$ and the twisted Brill-Noether loci. For $G(\alpha;n,d,k)$, the first four questions make sense, but Question 5 needs to be replaced by
\begin {qu} Is $G(\alpha;n,d,k)$ smooth?
\end{qu}

For $n=1$ (the classical case), the answers to the above questions are all known (for general $C$) (see \cite{acgh} for this and much more).
\begin{itemize}
\item[(i)] For any $C$, $G(1,d,k)$ and $B(1,d,k)$ are non-empty if $\beta(1,d,k)\ge0$ and connected if $\beta(1,d,k)>0$.
\item[(ii)] For general $C$, $G(1,d,k)$ is smooth of dimension $\beta(1,d,k)$ (in particular, it is empty if $\beta(1,d,k)<0$) and it is irreducible if $\beta(1,d,k)>0$.
\item[(iii)] For general $C$, $B(1,d,k)$ has dimension $\beta(1,d,k)$ if 
\[0\le\beta(1,d,k)\le g,\] is irreducible if $\beta(1,d,k)>0$, and its singular set coincides with $B(1,d,k+1)$ if $\beta(1,d,k)<g$.
\end{itemize}
Properties (ii) and (iii) follow from the fundamental fact that, for the general curve $C$ of genus $g$, the Petri map $\mu_E$ is injective for every line bundle $E$ on $C$. A curve with this property is said to be a \textit{Petri curve} (the terms \textit{Petri general} and \textit{Brill-Noether-Petri} are often used in the literature, \textit{Brill-Noether general} being reserved for a weaker concept).

For the twisted Brill-Noether loci $B(1,d,k)(F)$, all this remains true, provided that, in (ii) and (iii), we take $F$ to be general as well as $C$ \cite{gh,laz,te2,hi} (see also  \cite[Theorem 2.1]{hhn}).

If $n\ge2$, none of the above is true for all values of $(d,k)$ with the possible exception of connectedness. In fact, the following question remains (to my knowledge) open.
\begin{qu} Is it true that, if $B(n,d,k)\ne\emptyset$ and $\beta(n,d,k)>0$, then $B(n,d,k)$ is connected?
\end{qu}

There are some cases in which there are good and expected results for $B(n,d,k)$, for instance for $d\le2n$ \cite{bgn,m1,m2} (hence also $g\le3$) and hyperelliptic curves \cite{bmno}. There are already considerable gaps in our knowledge for non-hyperelliptic curves of genus $4$ \cite{ln1}. In general, the condition $\beta(n,d,k)\ge0$ is not sufficient for non-emptiness of $\beta(n,d,k)$ (see, for example, \cite{bgn}), nor is it necessary even when $C$ is general (examples can be obtained using the results of \cite{m3}). However there is a wide range of values for $(n,d,k)$ for which, if $C$ is general, there exists a component of $B(n,d,k)$ of the expected dimension and generically smooth; this was proved in \cite{te} (see \cite{clt} for the smoothness) using degenerations and in \cite{m3} using elementary transformations. There is an interesting question in connection with this.

\begin{qu} Are the components of $B(n,d,k)$ constructed in \cite{te} and \cite{m3} the same?
\end{qu}

There is a third construction in \cite{bmno}, which overlaps with those of \cite{te} and \cite{m3}. One can ask whether this also gives the same component when all are defined. 

 \begin{qu} What is a ``good'' definition for \textit{Brill-Noether general for rank $n$}?
 \end{qu}
 At first sight, this looks like a good question. However, since some of the answers to the basic questions are negative even on the general curve, and others are unknown, it is difficult to specify any specific conditions which must be satisfied on the general curve except by looking at specific values for $(d,k)$. This certainly allows us to state conditions which are necessary for a curve to be called Brill-Noether general, and one can also find sufficient conditions, but these are some distance apart.
 
 \begin{qu} What is a ``good'' definition of $\alpha$-stability for twisted coherent systems?
 \end{qu}
 This, on the other hand, is an excellent question. Recently, Schmitt has proposed two definitions for coherent systems of type $(\rho,n,d,k)(F)$, where $F$ is a line bundle. The obvious extension of the standard definition for standard coherent systems is one of these, but presents some problems over the permissible range for $\alpha$ \cite{sch1}. The second one is more complicated, but fits well with Schmitt's concept of decorated bundles and has some interesting features \cite{sch2}. A key requirement for a good definition is that it leads to a construction of moduli spaces.

\section{Clifford indices and gonalities}\label{ci}

In determining the possible dimensions for the spaces of sections of line bundles on $C$, the starting point is the classical theorems of Riemann-Roch and Clifford. In higher rank, the Riemann-Roch Theorem generalises for any vector bundle: if $E$ has rank $n$ and degree $d$, then 
\[h^0(E)-h^0(E^*\otimes K)=d-n(g-1).\]
For Clifford's Theorem, we have to restrict to semistable $E$, when we get, for $0\le d\le n(2g-2)$,
\[h^0(E)\le \frac{d}2+n.\]
Note that, if $E$ is semistable and $d<0$, $h^0(E)=0$, while, if $d>n(2g-2)$, $h^0(E)=d-n(g-1)$.

These results depend only on the genus of $C$. To take the complex structure of $C$ into account, we need more refined invariants. The first of these is the \textit{Clifford index}, defined as follows. For any vector bundle $E$ of rank $n$ and degree  $d$,
\[
\Cliff(E) := \frac{1}{n} (d - 2(h^0(E) -n)) = \mu(E) -2\frac{h^0(E)}{n} + 2.\]
For $g\ge4$, the \textit{Clifford indices} of $C$ are then defined by
\[\Cliff_n(C):= \min_{E} \left\{ \Cliff(E) \;\left| 
\begin{array}{c} E \;\mbox{semistable of rank}\; n \\
h^0(E) \geq 2n,\; \mu(E) \leq g-1
\end{array} \right\}. \right.
\]
Note that, using Serre duality, we can extend this definition to the full range $0\le d\le n(2g-2)$. For $n=1$, we get the classical Clifford index, which we denote here by $\Cliff_1(C)$. We can extend the definition to $g\ge2$ by defining $\Cliff_n(C)$ to be $0$ when $g=2$ and when $C$ is hyperelliptic of genus $3$, and to be $1$ when $C$ is non-hyperelliptic of genus $3$. (Clifford indices were defined in \cite{ln}, where they were denoted by $\gamma_n'$; the notation $\gamma_n$ has also been used, but this was used in \cite{ln} with a different meaning.)

We can define still more refined invariants 
\[\go_{n,k}:=\min\{d|\exists\text{ semistable }E\text{ of rank }n\text{ with }h^0(E)\ge k\}.\]
Classically, $\go_{1,2}$ is known as the \textit{gonality} of $C$ and we shall refer to all the $\go_{n,k}$ as \textit{gonalities}. We shall also define $\go_{n,k}^s$ to be the corresponding quantities when we replace ``semistable'' by ``stable'' in the definition.

Curves of genus $0$ or $1$ are excluded from the definition of Clifford index and all gonalities are known in these cases, so we assume for the rest of this section that $g\ge2$.

It is conventional to represent the non-emptiness of $B(n,d,k)$ graphically by plotting $\lambda:=\frac{k}n$ against $\mu$. On this diagram, the upper bound given by Clifford's Theorem is the line $\lambda=\frac{\mu}2+1$ (valid for $0\le\mu\le2g-2$). However, it is certainly possible to find smaller upper bounds. Mercat has proposed the following problems.

\begin{qu} Find a function $f$ (defined for $0\le\mu\le2g-2$) as small as possible such that all $\widetilde{B}(n,d,k)$ with $\lambda>f(\mu)$ are empty.
\end{qu}
\begin{qu}\label{qu11} Find a function $h$ (defined for $0\le\mu\le2g-2$) as large as possible such that all $B(n,d,k)$ with $\lambda<h(\mu)$ are non-empty.
\end{qu}

Precise answers are known for $g=2$, $g=3$ and hyperelliptic curves of genus $4$ \cite{ln1} and $5$ \cite{ln2}; in fact, in these cases, all non-empty Brill-Noether loci are known. For hyperelliptic curves of any genus, there are very good results; in fact, if one weakens the last question to insist only that $B(n,d,k)\ne\emptyset$ for some $(n,d,k)$ with $d=n\mu$ and $k=n\lambda$, then again there is a precise answer \cite{bmno}. For bielliptic curves, one can take $h(\mu)=\frac\mu2$ \cite{ba} and, over a certain part of the range for $\mu$, one can also take $f(\mu)=\frac\mu2$ \cite{m}. On any curve, in the range $0<\mu<2$, we have $f(\mu)=h(\mu)=\frac1g(\mu-1)+1$ \cite{bgn,m1}, with a corresponding result for $2g-4<\mu<2g-2$.

For fixed $n$, the best known general result for $f$ is $f(\mu)=\frac{\mu-\Cliff_n(C)}2+1$ given by the Clifford index, and, for $h$, the bound of Teixidor \cite{te} and Mercat \cite{m3} as extended in  \cite{bmno}. The bound for $f$ can certainly be improved (see, for example, the diagrams for $4\le g\le6$ in \cite{ln1,ln2,ln5} and the example in \cite{li} for $g=41$). All of these bounds are piecewise linear (but not necessarily continuous). However, for a large class of curves on K3 surfaces, a recent result \cite{fl} establishes the quadratic bound $f(\mu)<1+\frac{g}{4(g-1)^2}\mu^2+\frac1g$ for $\mu\le g-1$ and the positivity of the Brill-Noether number also takes a quadratic form (although it does not provide a bound). It is even possible that the ``true'' bounds are fractal.
 
For twisted Brill-Noether loci, non-emptiness  and smoothness results somewhat analogous to those of \cite{te} and \cite{m3} have been obtained in \cite{hhn}. However, there are extra numerical conditions required in the proof.
\begin{qu} Can the additional conditions in \cite[Theorem 1.2]{hhn} for the non-emptiness of $B(n,d,k)(F)$ and the existence of a component which is generically smooth of dimension $\beta(n,d,k)(F)$ be relaxed?
\end{qu}

We have clearly $\Cliff_n(C)\le\Cliff_1(C)$ and indeed $\Cliff_{am}(C)\le\Cliff_m(C)$ for any positive integer $a$. The following question seems natural.
\begin{qu} Is it true that $\Cliff_n(C)\le\Cliff_m(C)$ when $n\ge m$?  
\end{qu}

An implication of an affirmative answer to this question for a given $n$ is that every bundle computing $\Cliff_n(C)$ would then be primitive (i.e. both $E$ and $E^*\otimes K$ are generated) \cite{ln6}. However, the answer can be negative; for $n=3$, there are examples in \cite{ln9}.

In \cite{m}, Mercat made a conjecture proposing upper bounds for $h^0(E)$ in the range $1\le \mu(E)\le 2g-3$. The main implication of this conjecture can be stated as
\begin{conj}
$\Cliff_n(C)=\Cliff_1(C).$
\end{conj}
This has come to be known as \textit{Mercat's Conjecture}. The conjecture was originally made for general $C$, but it makes sense for any $C$.

\begin{qu} For which values of $n\ge2$ and for which curves does Mercat's Conjecture hold? More generally, find a good lower bound for $\Cliff_n(C)$.
\end{qu}

The conjecture always holds if $\Cliff_1(C)\le2$, i.e. for hyperelliptic curves, trigonal curves, tetragonal curves, including bielliptic curves, and smooth plane curves of degree $6$. It has been proved recently that it holds for arbitrary smooth plane curves \cite{fl}. For the class of curves on K3 surfaces mentioned above, the exact value of $\Cliff_n(C)$ is calculated for $g\ge n^2$ and, for $n\ge3$, it is strictly less than $\Cliff_1(C)$ except for $n=3$, $g=10,14$. When $n =2$, the conjecture also holds for  general $\ell$-gonal curves with $\ell\ge5$ and $g\ge4\ell-4$ (possibly for smaller $g$ as well) \cite{ln4}, for arbitrary $\ell$-gonal curves with $g\ge(\ell-1)(2\ell-4)$ \cite{fo2}, and for $C$ general \cite{baf}. However, there exist $C$ (even having maximal Clifford index $\Cliff_1(C)=\lfloor\frac{g-1}2\rfloor$) for which it does not hold \cite{fo,ln7,fo2}; there are indeed Petri curves for which $\Cliff_2(C)<\Cliff_1(C)$ \cite{lc}. Moreover, $\Cliff_2(C)\ge\min\{\Cliff_1(C), \frac{\go_{1,5}-4}2\}$ \cite{ln} and this bound can be attained for any value of $\Cliff_1(C)$ \cite{ln7}.   

For $n=3$, the conjecture fails for general $C$ of genus $9$ or $11$ \cite{lmn} and for any smooth curve lying on a K3 surface with $g=9$ or $g\ge11$ and $\Cliff_1(C)=\left\lfloor\frac{g-1}2\right\rfloor$ \cite{fo2}. For the best bound known to me which covers all curves, see \cite{ln9}. In  particular, if $\Cliff_1(C)=3$, we have $\frac83\le\Cliff_3(C)\le3$.

\begin{qu}\label{qu16} Do there exist curves with $\Cliff_1(C)=3$ and $\Cliff_3(C)=\frac83$?
\end{qu}

There are severe restrictions on the existence of such a curve \cite{ln16}; in particular, one must have $9\le g\le12$.

Turning now to looking at gonalities, we have the following general question.
\begin{qu} What are the values of $\go_{n,k}$, $\go^s_{n,k}$ for a general curve? When is $\go_{n,k}=\go^s_{n,k}$?
\end{qu}
All values of $\go_{1,k}=\go_{1,k}^s$ are known for Petri curves, and also for hyperelliptic curves, trigonal curves, general tetragonal curves, bielliptic curves and smooth plane curves. For $n\ge2$, the values of $d_{n,k}$ and $d_{n,k}^s$ are known for $g=2$ and $g=3$ and for hyperelliptic curves of genus $4$ or $5$. 

The gonality $\go_{n,2n}$ is particularly interesting as bundles of this gonality are the first candidates for computing the Clifford index. This happens if and only if $d_{n,2n}=n(\Cliff_n(C)+2)$. In fact, it is easy to see that, on a general curve of genus $g$, there exist strictly semistable bundles of any rank $n\ge2$ and degree $d= n(\Cliff_1(C)+2)=n\lfloor\frac{g+3}2\rfloor$, so this certainly holds if Mercat's Conjecture holds, in particular if $n=2$. On the other hand, it is not obvious that there exist stable bundles with this property. In general, one can ask the following question. 
\begin{qu} Which bundles compute the Clifford index? Which, if any, stable bundles do so?
\end{qu}
For hyperelliptic and trigonal curves, the bundles computing $\Cliff_n(C)$ are all strictly semistable. For hyperelliptic curves, they are all direct sums of copies of the hyperelliptic line bundle \cite{re}; for trigonal curves of genus $g\ge5$, they are direct sums of copies of the unique trigonal line bundle \cite{ln3,ho}. For $n=2$,\cite{ln4} contains several further results as follows. The only bundle computing $\Cliff_2(C)=\Cliff_1(C)$ for a smooth plane curve of degree $\delta\ge5$ is $H\oplus H$, where $H$ is the hyperplane bundle. For the general tetragonal curve of genus $g\ge27$, the only bundle computing $\Cliff_2(C)=\Cliff_1(C)$ is $Q\oplus Q$, where $Q$ is the unique tetragonal bundle. For $\ell\ge5$, a similar result holds for the general $\ell$-gonal curve of  genus $g>\max\{3\ell^2-8\ell+7,46\}$. These lower bounds can certainly be improved.

In rank $2$, for $C$ general, we have $d_{2,4}=d_{2,4}^s=g+3$ if $g$ is odd \cite{te9,fo} and $d_{2,4}=g+2$ if $g$ is even. This leaves the following question.
\begin{qu}\label{qu39} Does there exist a stable bundle of rank $2$ and degree $g+2$ with $h^0=4$ on a general curve of even  genus $g$?
\end{qu}
The answer is negative for $g\le10$, although there do exist Petri curves of genus $g=10$ for which $B(2,12,4)\ne\emptyset$ \cite{gmn}. To my knowledge, the answer to the question is not known for $g\ge12$.  To complete the picture for bundles of rank $2$ computing $\Cliff_2(C)$, we can ask
\begin{qu} Does there exist a bundle $E$ computing $\Cliff_2(C)$ on a general curve with $h^0(E)>4$?
\end{qu}
It is shown in \cite{baf} that there are no such bundles on a curve of even genus $g\ge10$, while, on a curve of odd genus $g\ge15$, the only possibility is for a semistable bundle $E$ of degree $2g-2$ with $h^0(E)=\frac{g+3}2$. It is not known whether any such bundle $E$ exists, but, if it does, $\det E\not\cong K$. For $6\le g\le9$, $g=11$ and $g=13$, there are bundles of determinant $K$ with the requisite number of sections. For $g=4$ and $g=5$, there are no candidates for such a bundle.

Much less is known when $n=3$.
 \begin{qu} What is the value of $\Cliff_3(C)$ for a general curve $C$ of genus $g$? Which bundles compute $\Cliff_3(C)$?
\end{qu}
Certainly $\Cliff_3(C)=\Cliff_1(C)$ for $g\le6$, but (to my knowledge) the answer is not known for $g=7$ (here $\frac83\le\Cliff_3(C)\le\Cliff_1(C)=3$). For $g=9$, we have $\Cliff_1(C)=4$ and $\Cliff_3(C)=\frac{10}3$ \cite{lmn,ln9}; moreover, in this case, $\Cliff_3(C)$ is computed by a stable bundle of degree $16$ with $h^0=6$. For $g=11$, it was shown in \cite{lmn,ln9} that $4\le\Cliff_3(C)<\Cliff_1(C)=5$; recently, the exact value $\Cliff_3(C)=\frac{14}3$ was determined in \cite{fl}. In general, the results of \cite{fl} imply that, for $C$ general of genus $g\ge9$, we have 
\[\frac23(g-1)-\frac23\left\lfloor\frac{g}3\right\rfloor\le\Cliff_3(C)\le\Cliff_1(C)=\left\lfloor\frac{g-1}2\right\rfloor.\]

There are many questions that can be asked for special curves. An interesting case is that of smooth plane curves of degree $\delta\ge7$. As seen above, we know that $\Cliff_n(C)=\Cliff_1(C)=\delta-4$, but, to my knowledge, the following question remains open.
\begin{qu} For $C$ a smooth plane curve of degree $\delta\ge7$, is it true that the only bundles computing $\Cliff_n(C)$ are direct sums of copies of the hyperplane bundle $H$?
\end{qu}

\section{Fixed determinant}\label{fd}
We again assume that $g\ge2$. Let $\cL$ be a line bundle of degree $d$. The na\um{\i}ve lower bound $\beta(n,d,k)-g$ for the dimension of a component of $G(\alpha;n,\cL,k)$ (and, when $k\ge d-n(g-1)$, that of a component of $B(n,\cL,k)$) cannot always be attained, although we do expect it to apply for the general $\cL$ of any fixed degree; for circumstances in which it can be attained in the case $n=2$, see \cite{te4}. Based on work of Osserman \cite{oss1,oss2}, we pose the following question.
\begin{qu}\label{qu41} Is it true that the dimension of $G(\alpha,n,\cL,k)$ at $(E,V)$ is at least $\beta(n,d,k)-g+h^1(\cL)\binom{k}{n}$? Can this lower bound be attained?
\end{qu}
Following Osserman's work, an affirmative answer to the first question has been given by Zhang \cite{zh2} for any Petri curve when $(E,V)$ is generically generated. Osserman also gave some examples in which the bound is attained and it is shown in \cite{gn} that the bound is attained for $k\le n+1$ under a mild generality condition on $\cL$. That some such condition is necessary is already implicit in Osserman's work. For $k>n+1$, the second question remains open.

For the rest of this section, we fix $n$ to be $2$. The case of $B(2,K,k)$ is particularly interesting. As indicated in Section \ref{bd}, there is a modified Brill-Noether number $\beta(2,K,k):=3g-3-\frac{k(k+1)}2$ in this case. Many years ago, Bertram and Feinberg \cite{bf} conjectured that, on a general curve $C$, $B(2,K,k)\ne\emptyset$ if and only if $\beta(2,K,k)\ge0$; the question was also asked by Mukai in the form of two problems \cite{mu1,mu2}. The ``only if'' part of this conjecture was proved by Teixidor \cite{te1}, who showed that the Petri map $\mu_E^K$ is injective for all semistable $E$ of rank $2$ and determinant $K$. It therefore remains to consider the following question.
\begin{qu} Is it true that $B(2,K,k)\ne\emptyset$ whenever $\beta(2,K,k)\ge0$?
\end{qu}
This question has been attacked by degeneration methods \cite{te3,zh1} and, more recently, by cohomological methods \cite{lnp}, which reduce the problem to a (complicated) combinatorial one. In fact, this shows that, if $B(2,K,k)$ is non-empty for one value of $g$, then this also holds for all greater values of $g$. The answer is now known to be affirmative for $k\le9$ and, correspondingly, for $g\le19$ (for the case $k=8$, $g=13$, see \cite{fjp}; $k=9$, $g=16$ is even more recent and, so far, unpublished). In the case $k=10$, we can reduce the problem to two particular values of $g$. In fact, $B(2,K,10)\ne\emptyset$ for $g\ge22$ (for $g=22$, see \cite{cgz}).
\begin{qu} Is $B(2,K,10)$ non-empty for $g=20$ or $g=21$?
\end{qu}

It is also known that $B(2,K,k)\ne\emptyset$ whenever $\beta(2,K,k)\ge0$ for $k=11$, $k=15$, $k=16$, $k=20$ and $k=24$; in fact, the calculations in \cite{lnp} and the papers mentioned above give complete answers for all genera $g<110$ except for $20,21,27,28,32-36,52,58,65,66,78,86-88,93-96$. The expectation is that these are not genuine exceptions.

A particularly interesting fact is that $\beta(2,K,k)$ can be greater than $\beta(2,2g-2,k)$, so there is the possibility that $B(2,K,k)$ is contained in a \textit{superabundant} component of $B(2,2g-2,k)$ (i.e., one of dimension greater than $\beta(2,2g-2,k)$). The first time this happens is for $g=5$  \cite{n}. This can be generalised to $B(2m,m(2g-2),k)$ by considering bundles of rank $2m$ with a $K$-valued symplectic structure $E\cong E^*\otimes K$ \cite{bh}; the required lower bound for the genus for this construction to give a superabundant component is  $g=50$.

We can also consider $B(2, \cL,k)$ when $\cL\not\cong K$. This case is of course covered by Question \ref{qu41}. It is possible to study $B(2,\cL,k)$ using cohomological methods similar to those of \cite{lnp}. This is done in the case when $d$ is odd in \cite{lns}. As for the case $\cL=K$, we reduce the problem (at least for general $\cL$) to a combinatorial one.

\section{Butler's Conjecture}\label{bu}

Suppose that $(E,V)$ is a generated coherent system on $C$ of type $(n,d,k)$ with $k>n$. We define a bundle $D_{E,V}$ by the exact sequence
\[0\lra D_{E,V}^*\lra V\otimes{\cO}_C\lra E\lra 0.\] 
Provided $h^0(E^*)=0$ (which is certainly the case if $E$ is semistable), dualising this sequence gives rise to a coherent system $D(E,V):=(D_{E,V},V^*)$. This is known as the \textit{dual span} construction. The bundle $D_{E,V}^*$ is often referred to as a \textit{kernel bundle} or \textit{syzygy bundle} and denoted by $M_{E,V}$. When $V=H^0(E)$, we write $D_E$ and $M_E$ for $D_{E,V}$ and $M_{E,V}$. Now let $S_0(n,d,k)$ be the open subset of $G_0(n,d,k)$ consisting of generated coherent systems. In \cite{bu2}, D. C. Butler made the following conjecture.

\begin{conj}\label{conj2}   Let $C$ be a general curve of genus $g$ and $n,d,k$ positive integers. Then, for a general $(E,V)\in S_{0}(n,d,k)$,
$D(E,V)\in S_{0}(k-n,d,k)$. Moreover,  $S_{0}(n,d,k)$ and $S_{0}(k-n,d,k)$ are birational.
\end{conj}
Here, by saying that $(E,V)$ is general, we mean that it belongs to some dense open subset of $S_0(n,d,k)$. Butler's Conjecture is often stated in the form of the following question.

\begin{qu} For a general $(E,V)\in S_0(n,d,k)$ on a general curve $C$ of genus $g\ge1$, is $D_{E,V}$ semistable?
\end{qu}
In general, this is slightly weaker than Conjecture \ref{conj2}, but is equivalent to it when $n=1$; note that, in this case, the hypothesis $(E,V)\in S_0(1,d,k)$ is equivalent to the simple assertion that $(E,V)$ is generated. This case has been examined by a number of authors and finally proved to be true in \cite{bbn2} . In fact, when $n=1$ and $g\ge3$, $D_{E,V}$ is very frequently stable, rather than just semistable. This is always true if $k\le5$ \cite{bbn1} and also if $k\ge6$ and $g\ge2k-6$. For $g=3$, the problem of stability would be completely solved by an answer to the following question.
\begin{qu} Suppose that $C$ is a non-hyperelliptic curve of genus $g=3$. Is it true that, for general $(E,V)\in G(1,2k,k)$ with $k\ge6$, $D_{E,V}$ is stable?
\end{qu}

An interesting concept, which dates back to Mumford, is that of \textit{linear stability} for a generated rank $1$ coherent system $(E,V)\in G(1,d,k)$. Linear stability is linked with the Chow stability of the image of $C$ in ${\mathbb P}^{k-1}$ under the map defined by $(E,V)$; for some recent work on this link, see \cite{bt}. Stability (semistability) of $D_{E,V}$ implies linear stability (semistability) of $(E,V)$, but the converse is not clear.

\begin{qu} Suppose that $(E,V)\in G(1,d,k)$ is generated and linearly stable (semistable). Is it true that $D_{E,V}$ is stable (semistable)?
\end{qu}

In \cite{ms}, this was shown to be true in many cases and some counterexamples were given. When $V=H^0(E)$, the question has been answered in the affirmative for Petri curves and hyperelliptic curves \cite{ct}, but counter-examples are known for smooth plane curves of genus $7$ \cite{cmt}. 

There is another condition on bundles, namely that of \textit{cohomological stability}. In fact, cohomological semistability is equivalent to semistability, but cohomological stability is stronger than stability. The concept has been used in proving Butler's Conjecture in some cases (see, for example, \cite{ms}).

Note that any generated coherent system of type $(n,d,n+1)$ is of the form $D(E,V)$ for some generated $(E,V)$ of type $(1,d,n+1)$. If one can prove that none of these $D(E,V)$ is stable, it follows that $B(n,d,n+1)$ contains no generated bundles. 

\begin{qu} Let $C$ be a general curve. Do there exist non-empty Brill-Noether loci $B(n,d,n+1)$ which contain no generated bundles?
\end{qu}

On special curves, such loci exist even for $n=1$. On a smooth plane curve of degree $7$, the example constructed in \cite{cmt} shows that $B(2,15,3)$, which is non-empty, contains no generated elements.  So far as I am aware, this particular example has not previously been observed.

For higher rank, much less is known. Conjecture \ref{conj2} is true on any curve of genus $g\ge2$ when $V=H^0(E)$ and $d>2ng$ (see \cite{bu1},\cite{m1}, \cite{bmgno}). It is also proved in \cite{bmgno} that the conjecture holds when $V=H^0(E)$ and $d=2ng$ provided $C$ is not hyperelliptic. Another interesting fact proved in \cite{bmgno} is that, if $(E,V)\in S_0(n,d,k)$ and  $D(E,V)\in S_0(k-n,d,k)$, then $(E,V)\in G(\alpha;n,d,k)$ for all $\alpha>0$, or, in other words, $(E,V)\in U^s(n,d,k)$; this is a possible source of counter-examples to Butler's Conjecture.

\begin{qu} Does Conjecture \ref{conj2} hold for $(E,V)\in S_0(2,d,4)$?
\end{qu}

This is the first case to consider for $n\ge2$ where $V\ne H^0(E)$. The conjecture  is shown to be true for a range of values of $d$ when $C$ is a general curve of genus $g\ge3$ in \cite{bmgno}, but there are many other values to be considered. 

\section{Coherent systems on the projective line}\label{p1}
It is well known that the only vector bundles on ${\mathbb P}^1$ are direct sums of line bundles ${\cO}_C(a)$. This means that the Brill-Noether theory is completely trivial. However, there is a very interesting theory of coherent systems (see \cite{ln11,ln12,ln13,nt}, also \cite{pp2}, where the closely related concept of holomorphic triples is studied).

The spaces $G(\alpha;n,d,k)$ are always smooth and irreducible of dimension $\beta(n,d,k)$ whenever they are non-empty. Moreover, the set $I(n,d,k):=\{\alpha|G(\alpha;n,d,k)\ne\emptyset\}$ is always an open interval (possibly semi-infinite) and all $G_i(\alpha;n,d,k)$ are birational.

\begin{qu} When is $I(n,d,k)$ non-empty? If $I(n,d,k)\ne\emptyset$ and we write $I(n,d,k)=]\alpha_m,\alpha_M[$, where $\alpha_M$ can be $\infty$, what are the values of $\alpha_m$ and $\alpha_M$?
\end{qu}

In many cases, the answer is known; for example, if $d=an$ for some integer $a\ge2$ and $k\ge n$, then $I(n,d,k)\ne\emptyset$ if and only if $\beta(n,d,k)=k((a+1)n-k)-n^2+1\ge0$. Moreover, if this condition holds, then $I(n,d,k)=]0,\infty[$. A complete answer is also known when $k\le3$. For $k<n$, there is a known upper bound for $\alpha_M$, so $\alpha_M\ne\infty$. For $k\ge n$, there are conjectures concerning the non-emptiness of $I(n,d,k)$ in \cite{nt} and, in particular, we can ask the following question.

\begin{qu} Suppose that $d=an-t$ with $a\ge2$, $1\le t\le n-1$ and $k\ge n$. If $G(\alpha;n,d,k)$ is non-empty for some $\alpha$, is it true that $\alpha_M=\infty$?
\end{qu}

I know of no cases in which the answer is no, but there is an indication for possible counter-examples in \cite{nt}.

\begin{qu}\label{qu63} Can one determine the flip loci for coherent systems of type $(n,d,k)$ on ${\mathbb P}^1$ and use this information to compute Hodge and Poincar\'e polynomials of $G_i(n,d,k)$?
\end{qu}

In the case $k=1$, this is done in \cite{ln13} and the information is sufficient to determine the Hodge polynomials of the smooth projective varieties $G_i(n,d,1)$. In the general case, the description of the flip loci in \cite{top1} works for the case $g=0$, but these loci are more complicated than in the case $k=1$ and it has not so far been possible to carry out the cohomological calculations.

\section{Coherent systems on elliptic curves}\label{ell}

When $g=1$, the basic questions for $G(\alpha;n,d,k)$ are fully answered in \cite{ln8} (see also \cite{pp1} for holomorphic triples). Precise conditions for non-emptiness are given and are entirely as expected, and all non-empty $G(\alpha;n,d,k)$ are smooth and irreducible of the expected dimension. As with ${\mathbb P}^1$, the set $I(n,d,k):=\{\alpha|G(\alpha;n,d,k)\ne\emptyset\}$ is always an open interval and all $G_i(\alpha;n,d,k)$ are birational. In fact, if it is non-empty,
\[I(n,d,k)=\begin{cases}]0,\frac{d}{n-k}[&\text{ if }k<n\\]0,\infty[&\text{ if }k\ge n.
\end{cases}\]  
There is a very interesting use of Fourier-Mukai transforms in \cite{ht}, which shows that $G_0(n,d,k)\cong G_0(n+ad,d,k)$ (and $G_L(n,d,k)\cong G_L(n+ad,d,k)$ if $k<n$) for all integers $a$.  In particular, the birational type of $G(\alpha;n,d,k)$ depends only on $n\bmod d$. The following question is raised in \cite{ln17}.

\begin{qu} Let $\cL$ be a line bundle of degree $d$. Is $G(\alpha;n,\cL,k)$ a rational variety?
\end{qu}

This is proved to be true when $\gcd(n,d)=1$ and in some other cases in \cite{ln17}.

The main question remaining is the analogue of Question \ref{qu63}.

\begin{qu}\label{qu71} When $g=1$, can one determine the flip loci for coherent systems of type $(n,d,k)$ and use this information to compute Hodge and Poincar\'e polynomials of $G_i(n,d,k)$?
\end{qu}

For some computations in this direction, covering also the fixed determinant case, see \cite{ln17}. In particular, the Hodge polynomial is computed for all $G_i(2+ad,d,1)$.

The results of \cite{ln8} can be used to obtain results on coherent systems on bielliptic curves \cite{ba2}.

\section{Coherent systems for $g\ge2$}\label{mcs}

We now assume that $g\ge2$. There is a version of Clifford's Theorem for $\alpha$-semistable coherent systems $(E,V)$ of type $(n,d,k)$ (just replace $h^0(E)$ by $k$ in the classical statement) \cite{ln14}. There are also obvious definitions of Clifford indices $\Cliff_{\alpha,n}(C)$, but I do not know of any work on the subject outside that of \cite{ln14}.

\begin{qu} What values can $\Cliff_{\alpha,n}(C)$ take on a curve with specified $\Cliff_1(C)$?
\end{qu}

Somewhat related to this question, we have

\begin{qu}
Find a good bound (either dependent on or independent of $\alpha$) for $h^0(E)$ when $(E,V)\in G(\alpha;n,d,k)$.
\end{qu}

It is known that such a bound exists \cite{top1}, but this bound can certainly be improved.

The major new feature for coherent systems is the variation of the moduli spaces with $\alpha$. Recall that $U^s(n,d,k)$ consists of those $(E,V)$ which are $\alpha$-stable for all $\alpha>0$. The following question was raised by Ballico in \cite{ba2}.

\begin{qu} What is the smallest integer $d'_{n,k}$ such that, for all $d\ge d'_{n,k}$, there exists a coherent system $(E,V)\in U(n,d,k)$ with $(E,V)$ generated and both the Petri map and the natural map $\psi:\bigwedge^n(V)\to H^0(\det E)$ injective?
\end{qu}

Ballico proved that $d'_{n,k}$ exists and gave a crude estimate for it. The significance of the injectivity of the natural map $\psi$ is that it is equivalent to saying that the image of $C$ under the map to projective space defined as the map to $\Gr(k,n)$ given by $(E,V)$  followed by the Pl\um{u}cker embedding is non-degenerate. It is clear that $d'_{n,k}\ge d_{n,k}$. The map $\psi$ is investigated in \cite{te6} in the case where $V=H^0(E)$ and injectivity is proved for general $E$ when $d=ng+1$ and $d=ng+2$. It follows that, in these cases, $\psi$ is injective for all $(E,V)$ with $V$ any subspace of $H^0(E)$. This suggests the possibility that $d'_{n,n+1}=ng+1$ and $d'_{n,n+2}=ng+2$.

A somewhat similar question (stated as a ``Model Theorem'') was proposed in \cite{n1}.

\begin{qu}\label{qu37} Suppose that $n\ge2$. Does there exist an integer $d_0(n,k)$ with the following properties
\begin{itemize}
\item[(a)]$G(\alpha;n,d,k)\ne\emptyset$ if and only if $\alpha>0$, $(n-k)\alpha<d$ and $d\ge d_0$;
\item[(b)]$B(n,d,k)\ne0$ if and only if $d\ge d_0$;
\item[(c)]if $d\ge d_0$, then $U(n,d,k)\ne\emptyset$?
\end{itemize}\end{qu}

For $k<n$, this is true and $d_0(n,k)=n-g(n-k)\ge0$ \cite{bgn,te5,top3}; moreover, if $d\le2n$, all the non-empty moduli spaces are irreducible of the expected dimension \cite{top3}. More generally, if we drop the ``only if'' from (a) and (b) and slightly modify (b), this is again true with the same value of $d_0(n,k)$ as is given by Teixidor's bound mentioned in connection with Question \ref{qu11} \cite{te5}.

For $k<n$, the moduli space $G_L(n,d,k)$ is described explicitly in \cite{bg} (see also \cite{top1}). There are also good results for $k=n$ and $k=n+1$ and some information for $k>n+1$. \cite{top1} contains also a description of the flip loci and estimates of their codimensions. If both flip loci have positive codimension, then there exist $(E,V)$ belonging to both $G_{i-1}(n,d,k)$ and $G_i(n,d,k)$. Some cases for low values of $k$ are worked out in detail. Very much related to this is the following question.

\begin{qu}Suppose that $C$ is general and $k>n$. Is it true that, if $G(\alpha;n,d,k)\ne\emptyset$ for some $\alpha>0$, then $G_L(n,d,k)\ne\emptyset$?
\end{qu}

One expects that this question usually has an affirmative answer, but there are counter-examples \cite{top3}. This shows in particular that the answer to Question \ref{qu37} may also be negative. Without any assumption about $C$, one can also ask the following question, which is simpler than the previous two questions.

\begin{qu} When is $U^s(n,d,k)$ $(U(n,d,k))$ non-empty?
\end{qu}

Papers which address this question directly include \cite{bp,bo,bo2,bmgno,bm}. When $C$ is general and $k=n+1$, the problem is related to Butler's Conjecture and is completely solved for $U^s(n,d,n+1)$ and solved in most cases for $U(n,d,n+1)$ \cite{bbn2}. For $n=2$, $k=4$, see \cite{te9,fo,bmgno}; here, $U(2,d,4)$ is non-empty for $g\ge4$ and $d\ge g+3$ and this is best possible for $g$ odd. For $g$ even, $U^s(2,g+2,4)=\emptyset$ (hence also $U(2,g+2,4)=\emptyset$) for $g\le10$. The following question (a version of Question \ref{qu39}) remains.

\begin{qu} For a general curve $C$ of even genus $g\ge12$, is $U^s(2,g+2,4)$ non-empty?
\end{qu}

Finally, we consider the question of a more detailed study of the wall-crossings, leading to a comparison of the cohomology of $G_{i-1}(n,d,k)$ and that of $G_i(n,d,k)$ and hopefully, by induction, relating the cohomology of $G_0(n,d,k)$ to that of $G_L(n,d,k)$.

\begin{qu} Can one obtain an expression for the change in the Hodge-Deligne polynomial (or Poincar\'e polynomial) when crossing the critical value $\alpha_i$?
\end{qu}

We use the term Hodge-Deligne polynomial because of the possibility of singularities in the moduli spaces. We have already discussed the case $g=0$ (Question \ref{qu63}, where there are complete results for $(n,d,1)$,  and $g=1$ (Question \ref{qu71}). For type $(2,\cL,1)$, with $\cL$ a line bundle of degree $d$, this has been worked out by Thaddeus \cite{th}, and it is easy to deduce the results for type $(2,d,1)$. In this case, the wall-crossings are genuine flips, which makes the calculations easier. These results were reproved in \cite{mov1}. The geometry of flips was discussed in the nicest case in \cite{top2} and results obtained for type $(n,d,n-2)$ with $n\ge3$ and critical values close to the upper bound $\frac{d}2$ for $\alpha$. For $(3,d,1)$, see \cite{mun}. The case $(4,d,1)$ is partially covered in the thesis of M. Tommasini and a start made on $(2,d,2)$ (unpublished - I hope that some of this work might be published soon). The related topic of triples of rank $(2,2)$ is discussed in \cite{mov2}.

\section{Final remarks}\label{fr}

In this section, we look at some topics not covered above.

\begin{rem}\label{r1}\begin{em}
One can study Brill-Noether loci for singular curves; the definition of coherent systems in \cite{kn} works for any polarised curve (connected scheme of pure dimension $1$ with a polarisation). Most work has concentrated on the case of a curve with nodal (or possibly cuspidal) singularities.  Brill-Noether theory on irreducible nodal and cuspidal curves has been studied, for example, in \cite{bh1,bhs}. In particular, the results of \cite{bgn} are generalised in \cite{bh1}, while kernel bundles are discussed in \cite{bhs}, including a proof of semistability (stability) in the case $d\ge2ng$ ($d>2ng$) (compare \cite{bu1}); \cite{bhs} also contains generalisations of the results of \cite{m1}, \cite{m2}. Coherent systems on integral curves are discussed in \cite{ba3,bap} and on irreducible nodal curves in \cite{bh2}, where results from \cite{bg} and \cite{top1} are generalised, including a description of $G_L(\alpha;n,d,k)$ when $k\le n$. The case of a nodal curve  of (arithmetic) genus $1$ is discussed in \cite{bh3}. Open questions are in general similar to those for smooth curves. One problem which has been solved in the smooth case, but remains open in the nodal case is the following.

\begin{qu} For a general $(E,V)\in G(1,d,k)$ on a general nodal curve, is $D_{E,V}$ semistable?
\end{qu}

The paper \cite{bhp} is interesting in its own right, but contains results which are necessary for solving this problem.

For a reducible nodal curve, the components of the moduli space $M(n,d)$ were determined in \cite{te7,te8} together with a useful criterion for the existence of stable bundles. For recent work on coherent systems on  reducible nodal curves, see \cite{brf2,brf3}. The failure of Butler's conjecture on a reducible nodal curve with $2$ components and one node is discussed in \cite{brf}.
\end{em}\end{rem}

\begin{rem}\label{r2}\begin{em}
The constructions in \cite{kn} also work in finite characteristic. There is no doubt that much of the theory applies in that case, but, to my knowledge, there has been no systematic work on this.
\end{em}\end{rem}

\begin{rem}\label{r3}\begin{em}
If $E$ is a bundle of rank $n$, one can define the Segre invariant of rank $r$ in terms of the maximal degree of a subbundle of rank $r$. There is substantial study of these invariants, including the stratifications of $M(n,d)$ which they determine \cite{bl,rt}. 
\begin{qu}\label{qu92} How are the Segre stratifications related to the stratification induced by the Brill-Noether loci in $M(n,d)$?
\end{qu}
Segre stratifications for the moduli space of coherent systems are defined in \cite{roa}. These should give information  on Question \ref{qu92} and can also be used to describe certain wall-crossings \cite{roa}.
\end{em}\end{rem}

\begin{rem}\begin{em}
The problem of Brill-Noether loci in rank $1$ for general curves of fixed gonality $m:=d_{1,2}$  has recently been solved \cite{llv}(following \cite{pfl,cpj,lar}). The methods involve degeneration to a chain of elliptic curves and some intricate combinatorial computations using the splitting type of the direct images of line bundles under the natural $m$-fold covering $C\to{\mathbb P}^1$.
\begin{qu}
Study the Brill-Noether loci and coherent systems in higher rank for a curve $C$ through the morphism $C\to{\mathbb P}^1$ induced by a line bundle of degree $d_{1,2}$.
\end{qu}

In view of the somewhat complicated arguments of \cite{llv}, this looks a hard problem in general. One could start with the trigonal case, where the rank-$1$ situation has been well understood for some time.
\end{em}\end{rem}

\begin{rem}\begin{em}
One can ask for Torelli theorems for Brill-Noether loci and moduli spaces of coherent systems, in other words, whether these loci determine the curve.
\begin{qu} To what extent do Brill-Noether loci and moduli spaces of coherent systems determine the curve $C$?
\end{qu}
\end{em}\end{rem}


\begin{thebibliography}{BGMNO00}
\bibitem[ACGH85]{acgh}
E. Arbarello, M. Cornalba, P. A. Griffiths and J. Harris: 
\emph{Geometry of Algebraic Curves I}. 
Grundlehren math. Wiss. 267, Springer, New York, 1985. 
\bibitem[BH21]{bh} A. Bajravani and G. H. Hitching: 
\emph{Brill-Noether loci on moduli spaces of symplectic bundles over curves}.
Collectanea Math. 72 (2021), 443--469.
\bibitem[BaF18]{baf} B. Bakker and G. Farkas:
\emph{The Mercat conjecture for stable rank $2$ vector bundles on generic curves}.
Amer. J. Math. 140 (2018), 1277--1295.
\bibitem[Ba98]{ba} E. Ballico:
\emph{Brill-Noether theory for vector bundles on projective curves}.
Proc. Camb. Phil. Soc. 124 (1998), 483--499.
\bibitem[Ba06a]{ba2} E. Ballico:
\emph{Coherent systems with many sections on projective curves}.
Internat. J. Math. 17 (2006), 263--267.
\bibitem[Ba06b]{ba3} E. Ballico:
\emph{Stable coherent systems on integral projective curves: an asymptotic existence theorem},
Int. J. Pure Appl. Math. 27 (2006), 205--214.
\bibitem[BaP07]{bap} E. Ballico and F. Prantil:
\emph{Coherent systems on singular genus one curves}.
Int. J. Contemp. Math. Sci. 2 (2007), 1527--1543, doi: 10.12988/ijcms.2007.07160.
\bibitem[BF98]{bf}
A.~Bertram and B.~Feinberg,
\emph{On stable rank two bundles with canonical determinant and many sections}.
In:  Algebraic Geometry (Catania, 1993/Barcelona 1994), 259--269, Lecture Notes in Pure and Appl. Math. {\bf 200},  Marcel Dekker, New York, 1998.
\bibitem[Bh07]{bh1} U. N. Bhosle: 
\emph{Brill-Noether theory on nodal curves}.
Internat. J. Math. 18 (2007), 1133--1150.
\bibitem[Bh09]{bh2} U. N. Bhosle:
\emph{Coherent systems on a nodal curve}.
In: Moduli spaces and vector bundles, London Mathematical Society Lecture Note Series Vol. 359, 437--455, Cambridge University Press, 2009.
\bibitem[Bh11]{bh3} U. N. Bhosle:
\emph{Coherent systems on a nodal curve of genus $1$}.
Math. Nachr. 284 (2011), 1829--1845, doi: 10.1002/mana.200910133.
\bibitem[BBN08]{bbn1} U. N. Bhosle, L. Brambila-Paz and P. E. Newstead:
\emph{On coherent systems of type $(n,d,n+1)$ on Petri curves}.
Manuscripta Math. 126 (2008), 409--441.
\bibitem[BBN15]{bbn2} U. N. Bhosle, L. Brambila-Paz and P. E. Newstead:
\emph{On linear series and a conjecture of D. C. Butler}.
Internat. J. Math. 26 (2015), 1550007 (18 pages).
\bibitem[BhP14]{bhp} U. N. Bhosle and A. J. Parameswaran:
\emph{Picard bundles and Brill- Noether loci on the compactified Jacobian of a nodal curve}.
 IMRN Vol. 2014, Issue 15 (2014), 4241--4290,
doi: 10.1093/imrn/rnt069.
\bibitem[BhS13]{bhs} U. N. Bhosle and S. K. Singh:
\emph{Brill-Noether loci and generated torsionfree sheaves over nodal and cuspidal curves}.
Manuscripta Math. 141 (2013), 241--271.
\bibitem[BG02]{bg} S. B. Bradlow and O. Garc\'{\i}a-Prada: 
\emph{An application of coherent systems to a Brill-Noether problem}.
J. Reine Angew. Math. 551 (2002), 123--143.
\bibitem[BGMN03]{top1} S. B. Bradlow, O. Garc\'{\i}a-Prada, V. Mu\~noz and P. E. Newstead: 
\emph{Coherent systems and Brill-Noether theory}. 
Internat. J. Math. 14 (2003), 683--733.
\bibitem[BGMMN07]{top2} S. B. Bradlow, O. Garc\'{\i}a-Prada, V. Mercat, V. Mu\~noz and P. E. Newstead: 
\emph{On the geometry of moduli spaces of coherent systems on algebraic curves}.
Internat. J. Math. 18 (2007), 411--453. 
\bibitem[BGMMN09]{top3} S. B. Bradlow, O. Garc\'{\i}a-Prada, V. Mercat, V. Mu\~noz and P. E. Newstead: 
\emph{Moduli spaces of coherent systems of small slope on algebraic curves}. 
Comm. in Alg. 37 (2009), 2649--2678.
\bibitem[BP08]{bp} L. Brambila-Paz:
\emph{Non-emptiness of moduli spaces of coherent systems}.
Internat. J. Math. 19 (2008), 779--799.
\bibitem[BGN97]{bgn} L. Brambila-Paz, I. Grzegorczyk and P. E. Newstead:
\emph{Geography of Brill-Noether loci for small slopes}.
J. Alg. Geom. 6 (1997), 645-669.
\bibitem[BL98]{bl} L. Brambila-Paz and H. Lange:
\emph{A stratification of the moduli space of vector bundles on curves (dedicated to Martin Kneser on the occasion of his 70th birthday)}.
J. Reine Angew. Math. 494 (1998), 173--187.
\bibitem[BM20]{bm} L. Brambila-Paz and O. Mata-Gutierrez:
\emph{$(t,\ell)$-stability and coherent systems}.
Glasgow Math. J. 62 (2020), 661--672.
\bibitem[BMGNO19]{bmgno} L. Brambila-Paz, O. Mata-Gutierrez, P. E. Newstead and A. Ortega:
\emph{Generated coherent systems and a conjecture of D. C. Butler}.
Internat. J. Math. 30 (2019), 1950024 (25 pages), doi: 10.1142/S0129167X19500241. 
\bibitem[BMNO00]{bmno} L. Brambila-Paz, V. Mercat, P. E. Newstead and F. Ongay:
\emph{Nonemptiness of Brill-Noether loci}.
Internat. J. Math. 11 (2000), 737--760.
\bibitem[BO09a]{bo} L. Brambila-Paz and A. Ortega:
\emph{Brill-Noether bundles and coherent systems on special curves}.
In: Moduli spaces and vector bundles, 456--472, London Mathematical Society Lecture Note Series Vol. 359,  Cambridge University Press, Cambridge, 2009.
\bibitem[BO09b]{bo2} L. Brambila-Paz and A. Ortega:
\emph{Estabilidad de sistemas coherentes}.
Aportaciones Mat. Comun. 40 (2009), 15--26.
\bibitem[BT16]{bt} L. Brambila-Paz and H. Torres-L\'opez:
\emph{On Chow stability for algebraic curves}.
Manuscripta Math. 151 (2016), 289--304.
\bibitem[BrF20a]{brf} S. Brivio and F. F. Favale:
\emph{On kernel bundles over reducible curves with a node}.
Internat. J. Math. 31 (2020), 2050054 (15 pages), doi: 10.1142/S0129167X12500371.
\bibitem[BrF20b]{brf2} S. Brivio and F. Favale:
\emph{Coherent systems on curves of compact type}.
J. Geom. Phys. 158 (2020), 103850, doi: 10.1016/j.geomphys.2020.103850.
\bibitem[BrF21]{brf3} S. Brivio and F. Favale:
\emph{Coherent systems and BGN extensions on nodal reducible curves}.
arXiv: 2104.06883
\bibitem[Bu94]{bu1} D. C. Butler:
\emph{Normal generation of vector bundles over a curve}.
J. Diff. Geom. 39 (1994), 1-34.
\bibitem[Bu97]{bu2} D. C. Butler:
\emph{Birational maps of moduli of Brill-Noether pairs}.
arXiv:alg-geom/9705009.
\bibitem[CLT18]{clt} A. Castorena, A. L\'opez Mart\'{\i}n and M. Teixidor i Bigas:
\emph{Petri map for vector bundles near good bundles}.
J. Pure Appl. Algebra. 222 (2018), 1692--1703.
\bibitem[CMT20]{cmt} A. Castorena, E. C. Mistretta and H. Torres-L\'{o}pez:
\emph{On linear stability and syzygy stability for rank $2$ linear series}.
arXiv: 2001.03609.
\bibitem[CT18]{ct} A. Castorena and H. Torres-L\'{o}pez:
\emph{Linear stability and stability of syzygy bundles}.
Internat. J. Math. 29 (2018), 1850080 (14 pages).
\bibitem[CPJ19]{cpj} K. Cook-Powell and D. Jensen:
\emph{Components of Brill-Noether loci for curves of fixed gonality}.
arXiv: 1907.08366.
\bibitem[CGZ21]{cgz} D. Cotterill, A. Alonso Gonzalo and N. Zhang:
\emph{The strong maximal rank conjecture and higher rank Brill-Noether theory}.
J. London Math. Soc. 104 (2021), 169--205.
\bibitem[FJP21]{fjp} G.Farkas, D. Jensen and S. Payne:
\emph{The non-abelian Brill-Noether divisor on $\overline{\cM}_{13}$ and the Kodaira dimension of $\overline{\cR}_{13}$.}
arXiv: 2110.09553.
\bibitem[FO11]{fo} G. Farkas and A. Ortega:
\emph{The maximal rank conjecture and rank two Brill-Noether theory}.
Pure and Appl. Math. Quarterly 7 no. 4 (2011) (Special Issue: In Memory of Eckart Viehweg), 1265--1296.
\bibitem[FO12]{fo2} G. Farkas and A. Ortega:
\emph{Higher rank Brill-Noether theory on sections of K3 surfaces}.
Internat. J. Math. 23 (2012), 1250075.
\bibitem[FP05]{fp} G. Farkas and M. Popa:
\emph{Effective divisors on $\overline{M}_g$, curves on K3 surfaces, and the slope conjecture}.
J. Alg. Geom. 14 (2005), 241--267.
\bibitem[FL21]{fl}  S. Feyzbakhsh and Chunyi Li: 
\emph{Higher rank Clifford indices of curves on a K3 surface}. 
Selecta Math. (N.S.)  27  (2021),  no. 3, Paper No. 48 (34 pages).
\bibitem[Ghi83]{gh} F. Ghione:
\emph{Un probl\`eme du type Brill-Noether pour les fibr\'es vectoriels}.
In: Algebraic Geometry - Open Problems (Ravello, 1982), 197--209, Lecture Notes in Mathematics Vol. 997, Springer, Berlin, 1983.
\bibitem[GMN11]{gmn} I. Grzegorczyk, V. Mercat and P. E. Newstead:
\emph{Stable bundles of rank 2 with four sections}.
 Internat. J. Math. 22 (2011), 1743--1762, doi: 10.1142/S0129167X11007434.
\bibitem[GN14]{gn} I. Grzegorczyk and P. E. Newstead:
\emph{On coherent systems with fixed determinant}.
Internat. J. Math. 25 (2014), 1450045 (11 pages), doi: 10.1142/S0129167X14500451. 
\bibitem[HT08]{ht} D. Hern\'andez Ruip\'erez and C. Tejero Prieto:
\emph{Fourier-Mukai transforms for coherent systems on elliptic curves}. 
J. London Math. Soc. 77 (2008), 15--32.
\bibitem[Hir88]{hi} A. Hirschowitz:
\emph{Probl\`emes de Brill-Noether en rang sup\'erieur}.
C. R. Acad. Sci. Paris S\'er. I Math. 307 (1988), 153--156.
\bibitem[HHN21]{hhn} G. H. Hitching, M. Hoff and P. E. Newstead: 
\emph{Nonemptiness and smoothness of twisted Brill-Noether loci}. 
Ann. Mat. Pura Applicata 200 (2021), 521--546, doi: 10.1007/s10231-020-01009-x (open access).
\bibitem[Hof21]{ho} M. Hoff:
\emph{A note on syzygies and normal generation for trigonal curves}.
arXiv: 2108.06106
\bibitem[KN95]{kn} A. King and P. E. Newstead:
\emph{Moduli of Brill-Noether pairs on algebraic curves}.
Internat. J.Math. 6 (1995), 733--748.
\bibitem[LMN12]{lmn} H. Lange, V. Mercat and P. E. Newstead: 
\emph{On an example of Mukai}. 
Glasgow Math. J. 54 (2012), 261--271, doi:10.1017/S0017089511000577.
\bibitem[LN04]{ln11} H. Lange and P. E. Newstead:
\emph{Coherent systems of genus $0$}.
Internat. J. Math. 15 (2004), 409--424.
\bibitem[LN05]{ln8} H. Lange and P. E. Newstead: 
\emph{Coherent systems on elliptic curves}.
Internat. J. Math. 16 (2005), 787--805.
\bibitem[LN07]{ln12} H. Lange and P. E. Newstead:
\emph{Coherent systems of genus $0$ II: Existence results for $k\ge3$}.
Internat. J. Math. 18 (2007), 363--393.
\bibitem[LN08a]{ln13} H. Lange and P. E. Newstead:
\emph{Coherent systems of genus $0$ III: Computation of flips for $k=1$}.
Internat. J. Math. 19 (2008), 1103--1119.
\bibitem[LN08b]{ln14} H. Lange and P. E. Newstead: 
\emph{Clifford's theorem for coherent systems}.  
Arch. Math. 90 (2008), 209--216, doi:10.1007/s00013-007-2534-3. 
\bibitem[LN09]{ln17} H. Lange and P. E. Newstead:
\emph{Hodge polynomials and birational types of moduli spaces of coherent systems on elliptic curves}. 
Manuscripta Math. 130 (2009), 1--19.
\bibitem[LN10a]{ln} H. Lange and P. E. Newstead: 
\emph{Clifford indices for vector bundles on curves}.
In: A. Schmitt (Ed.) Affine Flag Manifolds and Principal Bundles, 165--202, Trends in Mathematics, Birkh\"auser, Basel, 2010.
\bibitem[LN10b]{ln6} H. Lange and P. E. Newstead: 
\emph{Generation of vector bundles computing Clifford indices}.
Arch. Math. 94 (2010), 529--537, doi: 10.1007/S00013-017-0126-0.
\bibitem[LN11]{ln7} H. Lange and P. E. Newstead: 
\emph{Further examples of stable bundles of rank 2 with 4 sections}.
 Pure Appl. Math. Quarterly 7 no.4 (2011) (Special Issue: In Memory of Eckart Viehweg), 1517--1528.
\bibitem[LN13a]{ln3} H. Lange and P. E. Newstead:
\emph{Bundles computing Clifford indices on trigonal curves}.
Arch. Math. 101 (2013), 21--31, doi: 10.1007/S00013-013-0540-1.
\bibitem[LN13b]{ln4}H. Lange and P. E. Newstead:
\emph{Vector bundles of rank $2$ computing Clifford indices}.
Comm. in Algebra 41 (2013), 2317--2345.
\bibitem[LN13c]{ln9} H. Lange and P. E. Newstead:
\emph{On bundles of rank 3 computing Clifford indices}. 
Kyoto J. Math. 53, no. 1 (Memorial issue for Masaki Maruyama) (2013), 25--54.
\bibitem[LN13d]{ln16} H. Lange and P. E. Newstead:
\emph{Bundles of rank 3 on curves of Clifford index 3}. 
J. Symbolic Computation  57 (2013), 3--18, doi: 10.1016/j.jsc.2013.05.002.
\bibitem[LN17]{ln1}H. Lange and P. E. Newstead:
\emph{Higher rank BN-theory for curves of genus $4$}.
Comm. in Algebra 45 (2017), 3948--3966, doi: 10.1080/00927872.2016.1251938.
\bibitem[LN16]{ln2}H. Lange and P. E. Newstead:
\emph{Higher rank BN-theory for curves of genus $5$}.
Rev. Mat. Complut. 29 (2016), 691--717.
\bibitem[LN18]{ln5}H. Lange and P. E. Newstead:
\emph{Higher rank BN-theory for curves of genus $6$}. 
Internat. J. Math. 29 (2018), 1850014 (40 pages), doi: 10.1142/S0129167X18500143.
\bibitem[LNP16]{lnp} H. Lange, P. E. Newstead and Seong Suk Park.
\emph{Nonemptiness of Brill-Noether loci in M(2,K)}. 
Comm. in Algebra 44 (2016), 746--767, doi: 10.1080/00927872.2014.990020. 
\bibitem[LNS15]{lns} H.Lange, P. E. Newstead and V. Strehl:
\emph{Nonemptiness of Brill-Noether loci in M(2,L)}.
Internat. J. Math. 26 (2015), 1550108 (26 pages), doi: 10.1142/S0129167X15501086.
\bibitem[LLV20]{llv} E. Larson, H. K. Larson and I. Vogt:
\emph{Global Brill-Noether theory over the Hurwitz space}.
arXiv: 2008.10765.
\bibitem[Lar20]{lar} H. K. Larson:
\emph{Refined Brill-Noether theory for all trigonal curves}.
arXiv: 2002.00142.
\bibitem[Laz84]{laz} R. Lazarsfeld: 
\emph{Some applications of the theory of positive vector bundles}.
In: Complete intersections (Acireale, 1983), 29--61, Lecture Notes in Math. Vol. 1092, Springer, Berlin, 1984.
\bibitem[LC13]{lc} M. Lelli-Chiesa:
\emph{Stability of rank-$3$ Lazarsfeld-Mukai bundles on K3 surfaces}.
Proc. London Math. Soc. 107 (2013), 451--479, doi:10.1112/plms/pds087.
\bibitem[Li19]{li}  Chunyi Li: 
\emph{On stability conditions for the quintic threefold}. 
Invent. Math.  218  (2019), 301–-340.
\bibitem[Mer99a]{m1} V. Mercat: 
\emph{Le probl\`eme de Brill-Noether pour des fibr\'es stables de petite pente}.
J. reine angew. Math. 506 (1999), 1--41.
\bibitem[Mer99b]{m3} V. Mercat: 
\emph{Le probl\`eme de Brill-Noether et le th\'eor\`eme de Teixidor}.
Manuscripta Math. 98 (1999), 75--85.
\bibitem[Mer01]{m2} V. Mercat:
\emph{Fibr\'es stables de pente 2}.
Bull. London Math. Soc. 33 (2001), 535--542.
\bibitem[Mer02]{m} V. Mercat:
\emph{Clifford's theorem and higher rank vector bundles}.
Internat. J. Math. 13 (2002), 785--796.
\bibitem[MS12]{ms} E. C. Mistretta and L. Stoppino:
\emph{Linear series on curves: stability and Clifford index}.
Internat. J. Math. 23 (2012), 1250121 (25 pages).
\bibitem[Mu95]{mu1}
S. Mukai:
\emph{Vector bundles and Brill-Noether theory}. 
In Current topics in complex algebraic geometry (Berkeley, CA, 1992/93), 145--158, Math. Sci. Res. Inst. Publ., 28, Cambridge Univ.
Press, Cambridge, 1995.
\bibitem[Mu01]{mu2}
S. Mukai:
\emph{Non-abelian Brill-Noether theory and Fano 3-folds}. 
Sugaku Expositions 14 (2001), 125--153.
\bibitem[Mun08]{mun} V. Mu\~noz:
\emph{Hodge polynomials of the moduli spaces of rank $3$ pairs}.
Geom. Dedicata 136 (2008), 17--46.
\bibitem[MOV07]{mov1}  V. Mu\~noz, D. Ortega and M. J. V\'azquez-Gallo:
\emph{Hodge polynomials of the moduli space of pairs}.
Internat. J. Math. 18 (2007), 695--721.
\bibitem[MOV09]{mov2} V. Mu\~noz, D. Ortega and M. J. V\'azquez-Gallo:
\emph{Hodge polynomials of the moduli space of triples of rank $(2,2)$}.
Quart. J. Math. 60 (2009), 235--272.
\bibitem[New11]{n1} P. E. Newstead:
{Existence of $\alpha$-stable coherent systems on algebraic curves}. 
In: Grassmannians, Moduli Spaces and Vector Bundles, 121--139, Clay Math. Proc. Vol. 14, American Mathematical Society, Providence, RI, 2011.
\bibitem[New18]{n} P. E. Newstead:
\emph{Some examples of rank-2 Brill-Noether loci}.
Rev. Mat. Complut. 21 (2018), 201--215, doi: 10.1007/s13163-017-0241-6.
\bibitem[NT21]{nt} P. E. Newstead and M. Teixidor i Bigas:
\emph{Coherent systems on the projective line}.
Quart. J. Math. 72 (2021), 115--136, doi:10.1093/qmathj/haaa024.
\bibitem[Oss13a]{oss1} B. Osserman:
\emph{Brill-Noether loci with fixed determinant in rank $2$}.
Internat. J. Math. 24 (2013), 1350019 (24 pages).
\bibitem[Oss13b]{oss2} B. Osserman:
\emph{Special determinants in higher-rank Brill-Noether theory}.
Internat. J. Math. 24 (2013), 1350084 (20 pages).
\bibitem[PP07]{pp1} S. Pasotti and F. Prantil:
\emph{Holomorphic triples on elliptic curves}.
Results Math. 50 (2007), 227--239.
\bibitem[PP08]{pp2} S. Pasotti and F. Prantil:
\emph{Holomorphic triples of genus $0$}.
Central European J. Math. 6 (2008),  129--142.
\bibitem[Pfl17]{pfl} N. Pflueger:
\emph{Brill-Noether varieties of $k$-gonal curves}.
Adv. Math. 312 (2017), 46--63.
\bibitem[Re98]{re} R. Re:
\emph{Multiplication of sections and Clifford bounds for stable vector bundles on curves}.
Comm. in Alg. 26 (1998), 1931-1944. 
\bibitem[Roa20]{roa} L. Roa-Leguizam\'on:
\emph{Segre invariant and a stratification of the moduli space of coherent systems}.
Internat. J. Math. 31 (2020), 2050117 (32 pages).
\bibitem[RT99]{rt} B. Russo and M. Teixidor i Bigas:
\emph{On a conjecture of Lange}.
J. Alg. Geom. 8 (1999), 483--496.
\bibitem[Sch19]{sch1} A. Schmitt:
\emph{A general notion of coherent systems}.
hal-02391836 (39 pages).
\bibitem[Sch21]{sch2} A. Schmitt:
\emph{Notes on coherent systems}.
Rev. Mat.: Teoria y Apl. 28 (2021), 1--38.
\bibitem[Tei91a]{te} M. Teixidor i Bigas:
\emph{Brill-Noether theory for stable vector bundles}.
Duke Math. J. 62 (1991), 385--400.
\bibitem[Tei91b]{te7} M. Teixidor i Bigas:
\emph{Moduli spaces of (semi)stable vector bundles on tree-like curves}.
Math. Ann. 290 (1991), 341--348.
\bibitem[Tei95]{te8} M. Teixidor i Bigas:
\emph{Moduli spaces of vector bundles on reducible curves}.
Amer. J. Math. 117 (1995), 125--139.
\bibitem[Tei98]{te6} M. Teixidor i Bigas:
\emph{Curves in Grassmannians}.
Proc. Amer. Math. Soc. 126 (1998), 1597--1603.
\bibitem[Tei04]{te3} M. Teixidor i Bigas:
\emph{Rank two vector bundles with canonical determinant}.
Math. Nachr. 265 (2004), 100--106.
\bibitem[Tei07]{te9} M. Teixidor i Bigas:
\emph{Existence of coherent systems of rank two and dimension four}. 
Collect. Math. 58 (2007), 193--198.
\bibitem[Tei08a]{te1} M. Teixidor i Bigas:
\emph{Petri map for rank two bundles with canonical determinant}.
Compos. Math. 144 (2008), 705--720.
\bibitem[Tei08b]{te5} M. Teixidor i Bigas:
\emph{Existence of coherent sysytems II}.
Internat. J. Math. 19 (2008), 1269--1283.
\bibitem[Tei10]{te4} M. Teixidor i Bigas:
\emph{Existence of vector bundles of rank two with fixed determinant and sections}.
Proc. Japan Acad. 86 Ser. A (2010), 113--118.
\bibitem[Tei14]{te2} M. Teixidor i Bigas: 
\emph{Injectivity of the Petri map for twisted Brill-Noether loci}.
Manuscripta Math. 145 (2014), 389--397.
\bibitem[Th94]{th} M. Thaddeus:
\emph{Stable pairs, linear systems and the Verlinde formula}.
Invent. Math. 117 (1994), 317--353.
\bibitem[Zh16]{zh1} N. Zhang:
\emph{Towards the Bertram-Feinberg-Mukai Conjecture}.
J. Pure Appl. Algebra 220 (2016), 1588--1654.
\bibitem[Zh17]{zh2} N. Zhang:
\emph{Expected dimensions of higher rank Brill-Noether loci}.
Proc. Amer. Math. Soc. 145 (2017), 3735--3746.
\end{thebibliography}
\end{document}